\documentclass{article}
\usepackage[utf8]{inputenc}
\usepackage{lastpage}
\usepackage[english]{babel}
\usepackage{indentfirst}	
\usepackage{color}
\usepackage{graphicx}
\usepackage{tikz-cd}
\usepackage{microtype} 
\usepackage{hyperref}
\usepackage{xurl}
\usepackage{latexsym}
\usepackage{amssymb}
\usepackage{amsmath}
\usepackage{amsthm}
\usepackage{tikz}
\usepackage{tikz-cd}
\newtheorem{thm}{Theorem}[section]
 \newtheorem{cor}{Corollary}[section]
 
 \newtheorem{prop}{Proposition}[section]
 
\newtheorem{rem}{Remark}[section]
 \newtheorem{ex}{Example}[section]
 
 \newtheorem{quest}{Question}[section]
 \newtheorem{cla}{Claim}[section]
 \usepackage[num]{abntex2cite}

\providecommand{\keywords}[1]
{
  \small	
  \textbf{\textit{Keywords---}} #1
}
 
\citebrackets[]
\title{Absorbing ideals: a survey on $\omega$-stable groups}

\author{Bruno Moreira Fernandes}

\begin{document}

\maketitle
\begin{abstract}
    We introduce a new concept in the \textit{Absorbing Ideal Theory} in commutative rings, that is, the $\omega$-stable groups. We will provide examples and non-examples of these groups, and establish their relationship with $H$-congruence. Ultimately, we will study the action of these groups considering the transitivity condition.
\end{abstract}

\begin{center}

.
\end{center}

\section{Introduction}
Throughout this paper, we write $R$ to be commutative ring with identity, \texttt{Prp}$(R)$ to be the set of all proper ideals of $R$, \texttt{Max}$(R)$ to be the set of all maximal ideals of $R$, \texttt{Rd}$(R)$ to be the set of all radical ideals of $R$ and $\texttt{Aut}(R)$ the group of automorphisms of $R$. In order to facilitate the development of this work, we write $\mathbb{N}$ to denote the set of all positive integers.

 Absorbing ideals in commutative rings were studied, systematically, by D. Anderson and A. Badawi in \cite{and-bad}. They studied the notion of ``how far an ideal is from being prime" through a function $\omega_{R}$ (we only write $\omega$ when there is no doubt about which ring is in question) that associates each proper ideal $I$ of $R$ to the smallest positive integer for which $I$ is an absorbing ideal, something we will define precisely later.
 
 One of the results obtained in this theory that interests us is the following:
 \\
 \begin{cla} (\cite{and-bad}, Theorem 4.2, item (b))
  Let $f:R\longrightarrow T$ be a surjective ring homomorphism and $I$ a $n$-absorbing ideal of $R$ containing $ker(f)$. Then $f(I)$ is an $n$-absorbing ideal of $T$ if and only if $I$ is an $n$-absorbing ideal of $R$. Moreover, $\omega_{T}(f(I))=\omega_{R}(I)$. In particular, this holds if $f$ is an isomorphism.
\end{cla}

We consider this theorem as the most natural starting point for investigating subgroups of permutations of sets of ideals that are ``$\omega$-stable", concept that we will define later. In this work, we define this notion and provide some examples and non-examples of this property in certain rings. 

In the follow section, we present the definitions and basic notions about the \textit{Absorbing Ideal Theory} in commutative rings, some basic results of this theory and also some theoretical constructions.

After, we define the $\omega$-stable groups over sets of ideals of a ring $R$, providing some examples and non-examples of these groups in different rings. We will also define the $\omega$-congruence and establish their relationship to $\omega$-stable groups. Indeed, the existence of a $H$-congruence is a necessary condition for the existence of a $\omega$-stable group (see the Corollary \ref{cor3.2}).

Finally, we will study $\omega$-stable groups that act transitively. We will show that, for integral domains that are not fields, it is not necessary to consider the transitive property when wanting to investigate $\omega$-stable groups on the set of all proper ideals of the ring (see the Corollary \ref{cor4.3}).

\section{Basic definitions and some ring-theoretic constructions with absorbing ideals}

Let $R$ be a ring. We say that a proper ideal $I$ of $R$ is $n$-\textbf{absorbing} if, for all $x_1,\ldots,x_{n+1}\in R$ such that $x_1 \cdots x_{n+1}\in I$, there are $1\leq i_1<\cdots< i_n\leq n$ such that $x_{i_1}\cdots x_{i_n}\in I$. In this case, we say that $I$ is, simply, an \textbf{absorbing ideal}.

The most trivial example is a proper ideal $P\subset R$ that is a $1$-absorbing if and only if is prime.

Next, we present the main properties of absorbing ideals.
\\
\begin{prop} (\cite{and-bad}, Theorem 2.1).
   If $R$ is a ring, then the following statements are valid:
   \\
   \\
   $(a)$\,\,$I\subset R$ is an $n$-absorbing ideal if and only if, for every integer $m$, with $m>n$, and all $x_1,\ldots ,x_m\in R$ such that $x_1\cdots x_m \in I$, there are $1\leq i_1<\cdots <i_n\leq m$ such that $x_{i_1}\cdots x_{i_n} \in I$.
   \\
   \\
   $(b)$\,\,If $I\subset R$ is an $n$-absorbing ideal, then $I$ is $m$-absorbing for all $m\geq n$.
   \\
   \\
   $(c)$\,\,If $I_j$ is an $n_j$-absorbing ideal of $R$ for all $1\leq j\leq m$, then $I_1\cap \cdots \cap I_m$ is an $n$-absorbing ideal of $R$, where $n=n_1+\cdots +n_m$. In particular, if $P_1,\ldots, P_n$ are prime ideals of $R$, then $P_1\cap \cdots \cap P_n$ is a $n$-absorbing ideal of $R$.
   \\
   \\
   $(d)$\,\,If $R$ is an integral domain and $p_1,\ldots,p_n$ are prime elements of $R$, then $I=\left\langle p_1\cdots p_n\right \rangle$ is an $n$-absorbing ideal of $R$.
   \\
   \\
   $(e)$\,\,If $I\subset R$ is an $n$-absorbing ideal, then $\sqrt{I}$ is an $n$-absorbing ideal of $R$. Furthermore, $x^n\in I$, for all $x\in I$.
\end{prop}

To increase our power of investigation into absorbing ideals, we need to introduce a function that measures ``\textit{how far a proper ideal is from being prime}". Let $I$ be a proper ideal of $R$. We write $\omega_{R}(I)$ for denote the smallest positive integer for which $I$ is an absorbing ideal. When $\omega_{R}(I)$ does not exist, we will write $\omega_{R}(I)=\infty$; hence
 $I$ is not $n$-absorbing for every $n\in\mathbb{N}$. If $\omega_{R}(I)<\infty$, $I$ is $\omega_{R}(I)$-absorbing, but it is not $k$-absorbing for every $1\leq k<\omega_{R}(I)$. We can establish a well defined correspondence $\omega_{R}:\texttt{Prp}(R)\longrightarrow \mathbb{N}\cup \lbrace \infty\rbrace$ that associates each proper ideal $I\subset R$ to $\omega_{R}(I)$. To make things easier, we will write ``$\omega$" instead of ``$\omega_{R}$" whenever there is no ambiguity about which ring is in question.

 Below we present some examples of absorbing and non absorbing ideals.
\begin{ex}
    A proper ideal $P$ of $R$ is prime if and only if $\omega(P)=1$.
\end{ex}
\begin{ex}\label{ex2.2}
    In $\mathbb{Z}$, every ideal of the form $I_n=\left\langle p_1\cdots p_n\right\rangle$, with $p_1,\ldots, p_n \in \mathbb{Z}$ prime number, is an $n$-absorbing ideal. Furthermore, $\omega(I_n)=n$. Conversely, if $I$ is a proper ideal of $\mathbb{Z}$, then $I=I_n$, for some collection of prime numbers $p_1,\ldots, p_n\in \mathbb{Z}$ and, therefore , $\omega(I)=n$.
\end{ex}
\begin{ex}
    Recall that $\mathbb{Z}_{n}$ is an integral domain iff $n$ is prime. Hence, for example, $\mathbb{Z}_{4}$ is not an integral domain. Moreover, $\omega(\mathbf{0})=\infty$. Indeed, $\Bar{2}^{k}\in \mathbf{0}$, for every $k\geq 2$, but $\Bar{2}\not\in \mathbf{0}$.
\end{ex}

Let $\mathcal{I}$  be a non-empty set of proper ideals of $R$ and
$$
 \Omega(R,\mathcal{I}):=\lbrace \omega(I);\,\,I\in \mathcal{I}\rbrace.
 $$
 \begin{rem}
     Remark that $\Omega(R,\mathcal{I})=\omega(\mathcal{I})$. If $\mathcal{I}=$ \texttt{Prp}$(R)$, we write $ \Omega(R,\mathcal{I})=\Omega(R)$. Hence, we have $\Omega(R,\mathcal{I})\subseteq\Omega(R)$ always. Futhermore, $1\in \Omega(R)$ always, since every proper ideal of $R$ is contained in some maximal ideal (this is basically obtained by applying Zorn's Lemma to the set \texttt{Prp}$(R)$, which is partially ordered by considering the inclusion relation).
 \end{rem}

Follow, a classification of $\Omega(R)$ for some classes of rings $R$.
\\
\begin{thm} (\cite{and-bad}, Theorem 5.11).
    Let $R$ be a ring. Then
    \\
    \\
    $(i)$\,\,If $\texttt{Max}(R) = n <\infty$, then $\lbrace 1,\ldots,n\rbrace\subseteq \Omega(R)$. If $\texttt{Max}(R)$ is infinite, then $\mathbb{N}\subseteq \Omega(R)$;\\
    \\
    $(ii)$\,\,If $R$ is an artinian ring, then $\Omega(R)=\lbrace 1,\ldots, n\rbrace$, for some positive integer $n$. In particular, $\Omega(R)=\lbrace 1\rbrace$ iff $R$ is a field;\\
    \\
    $(iii)$\,\,If $R$ is a noetherian ring with Krull dimension $\geq 1$, then $\Omega(R)=\mathbb{N}$.\\
    \\
    $(iv)$\,\,Let $R$ be a valuation domain (not a field). Then $\Omega(R)=\mathbb{N}$ if $R$ is a discrete valuation domain (DVR). If $R$ is not a DVR, then $\Omega(R)=\lbrace 1,\infty\rbrace$ if all nonzero prime ideals of $R$ are idempotent and $\Omega(R)=\mathbb{N}\cup\lbrace \infty\rbrace$ if $R$ has a non-idempotent nonzero prime ideal.
\end{thm}

Let $I,J\in \mathcal{I}$. We define the relation $\sim_{\omega}$ in $\mathcal{I}$ by
$$
I\sim_{\omega} J \Longleftrightarrow \omega(I)=\omega(J)
$$
The relation $\sim_{\omega}$ on $\mathcal{I}$ is an equivalence relation, since the equality relation $=$ it is also. Hence, two ideals $I$ and $J$ of $R$ are $\omega$-\textbf{equivalent} when $\omega(I)=\omega(J)$.

For each ideal $I\in\mathcal{I}$, let $[I]_{\omega}$ be its equivalence class corresponding to this relation which, here, will be called $\omega$-\textbf{orbit} of $I$. Denote by $\dfrac{\mathcal{I}}{\sim_{\omega}}$ the set of all $\omega$-orbits of ideals of $\mathcal{I}$. We call this set $\dfrac{\mathcal{I}}{\sim_{\omega}}$ a $\omega$-\textbf{orbit set}.

\begin{rem}\label{obs2.2}
 There exist an one-to-one correspondence between $\Omega(R,\mathcal{I})$ and $\dfrac{\mathcal{I}}{\sim_{\omega}}$.
\end{rem}

\section{On $\omega$-stable groups}
The following result shows that there is a certain type of ``stability" in relation to absorbing ideals.
\\
\begin{thm}\label{teo3.1} Let $f:R\longrightarrow T$ be a surjective ring homomorphism and $I$ be an absorbing ideal of $R$ containing $ker(f)$. Then $f(I)$ is an $n$-absorbing ideal of $T$ if and only if $I$ is an $n$-absorbing ideal of $R$. Moreover, $\omega_{T}(f(I))=\omega_{R}(I)$. In particular, this holds if $f$ is an isomorphism.
\begin{proof}
    Let $I$ be an absorbing ideal of $R$ containing $ker(f)$. Suppose that $f(I)$ is an $n$-absorbing ideal of $T$. Let $x_1,\ldots,x_{n+1}\in R$ such that $x_1\cdots x_{n+1}\in I$. Then $f(x_1),\ldots,f(x_{n+1})\in T$ and $f(x_1)\cdots f(x_{n+1})\in f(I)$. Since $f(I)$ is an $n$-absorbing ideal, then there are $1\leq i_1<\cdots <i_{n}\leq n+1$ such that $f(x_{i_1})\cdots f(x_{i_n})\in f(I)$. Hence, there exist $i\in I$ such that $f(x_{i_1})\cdots f(x_{i_n})=f(i)$. Since $f$ is a ring homomorphism, then $x_{i_1}\cdots x_{i_n}-i\in ker(f)\subseteq I$ and hence, $x_{i_1}\cdots x_{i_n}\in I$. Thus $I$ is an $n$-absorbing ideal. 
    
    Reciprocally, suppose that $I$ is an $n$-absorbing ideal of $R$. Let $y_1,\ldots,y_{n+1}\in T$ such that $y_1\cdots y_{n+1}\in f(I)$. Then there exist $i\in I$ such that $y_1\cdots y_{n+1}=f(i)$. Since $f$ is surjective, there are $x_1,\ldots, x_{n+1}\in R$ such that $y_i=f(x_i)$ for every $i=1,\ldots, n+1$. Hence, we obtain $f(x_1)\cdots f(x_{n+1})=f(i)$. Since $f$ is a ring homomorphism, then $x_{i_1}\cdots x_{i_n}-i\in ker(f)\subseteq I$ and hence, $x_{i_1}\cdots x_{i_n}\in I$. Being $I$ an $n$-absorbing ideal, there are $1\leq i_1<\cdots <i_{n}\leq n+1$ such that $x_{i_1}\cdots x_{i_n}\in I$. So, $y_{i_1}\cdots y_{i_n}=f(x_{i_1})\cdots f(x_{i_n})=f(x_{i_1}\cdots x_{i_n})\in f(I)$. Therefore, $f(I)$ is an $n$-absorbing ideal of $T$.

    The ``moreover" statement is true because if it were $\omega_{T}(f(I))\neq \omega_R(I)$, then we would obtain a contradiction with what we have just proved.

    The ``particular" statement is obvious.
    
\end{proof}
\end{thm}
Given a subgroup $G\subseteq \texttt{Aut}(R)$, let $\mathcal{I}$ be a $G$-\textbf{invariant} set of absorbing ideals of $R$, that is, $g(I)\in \mathcal{I}$, for all $g\in G$ and all ideals $I\in \mathcal{I}$.
In this case, each $g\in G$ induces a bijection $\Tilde{g}:\mathcal{I}\ni I\longmapsto \Tilde{g}(I):=g(I)\in \mathcal{I}$. Thus, with $\texttt{Perm}(\mathcal{I})$ being the group of permutations of $\mathcal{I}$, we see that the set

$$
H_{G}(\mathcal{I}):=\lbrace h\in \texttt{Perm}(\mathcal{I})\,|\,h=\Tilde{g},\,\textrm{for some}\,g \in G\rbrace
$$
is, of course, a subgroup of $\texttt{Perm}(\mathcal{I})$. Furthermore, there is a natural left group action of $H_{G}(\mathcal{I})$ over $\mathcal{I}$. We denote by $\dfrac{\mathcal{I}}{H_{G}(\mathcal{I})}$ the set of $H_{G}(\mathcal{I})$-orbits
$$
\texttt{Orb}_{H_{G}(\mathcal{I})}(I):=\lbrace h(I)\,|\,h\in H_{G}(\mathcal{I})\rbrace
$$
with $I\in\mathcal{I}$.

Note that the group $H_G(\mathcal{I})$ ``acts stably" on the $\omega$-orbits of $\mathcal{I}$. In fact, for any $h\in H_{G}(\mathcal{I})$ and $\omega$-orbit $[I]_{\omega}$, we have
$$
J\in [I]_{\omega} \Leftrightarrow \omega(J)=\omega(I) \Leftrightarrow \omega(h(J))\stackrel{\ast}{=}\omega(J)=\omega(I) \Leftrightarrow h(J)\in [I]_{\omega}.
$$
\begin{rem}
    The equality $\ast$ follows from the fact that every $h\in H_{G}(\mathcal{I})$ is uniquely defined by an automorphism $g\in G$ of $R$, that is, $h=\tilde{g}$ for some automorphism $g\in G$, and making use of the particular statement of Theorem \ref{teo3.1}.
\end{rem}
 Note that all $H_{G}(\mathcal{I})$-orbits of $\mathcal{I}$ are contained in the $\omega$-orbits of $\mathcal{I}$, that is, $ \texttt{Orb}_{H_{G}(\mathcal{I})}(I)\subseteq [I]_{\omega}$, for all $I\in \mathcal{I}$.

More generally, if $\mathcal{I}$ be a non-empty set of proper ideals of $R$, we say that a subgroup $H$ of $\texttt{Perm}(\mathcal{I})$ is $\omega$-\textbf{stable} over $\mathcal{I}$, or simply $\omega$-\textbf{stable}, if every $H$-orbit $\texttt{Orb}_{H}(I)$ is contained in the $\omega$-orbit $[I]_{\omega}$, for all $I\in \mathcal{I}$. Under these conditions, it is immediate that $H$ is $\omega$-stable over $\mathcal{I}$ iff $\omega(h(I))=\omega(I)$, for all $h\in H$ and all $I\in \mathcal{I}$.

\begin{rem}\label{obs3.2}
    According to this definition, if $\mathcal{I}$ is a non-empty $G$-invariant set of absorbing ideals of $R$, for some subgroup $G\subseteq \texttt{Aut}(R)$, then $H_{G}(\mathcal{I})$ is $\omega_R$-stable over $\mathcal{I}$. This tells us that the $\omega$-orbits of $\mathcal{I}$ are $H_{G}$-invariant. Therefore, there is a natural left group action of $H_{G}(\mathcal{I})$ over the $\omega$-orbits of $\mathcal{I}$. 
    
    More generally, let $H$ is a subgroup of $\texttt{Perm}(\mathcal{I})$ of some set $\mathcal{I}$ of proper ideals of $R$. If $H$ is $\omega$-stable on $\mathcal{I}$, then there is a natural left group action of $H$ over the $\omega$-orbits of $\mathcal{I}$.
    
    Furthermore, it is clear that if $H$ is the trivial group, then $H$ is $\omega$-stable over any non-trivial set of proper ideals of $R$ and if $\infty\not\in\Omega(R)$, then $H_{G}(\texttt{Prp}(R))$ is $\omega$-stable over $\texttt{Prp}(R)$, for every group of automorphisms $G$ of $R$.
\end{rem}

\begin{ex}\label{ex3.1}
Recall that an automorphism $\sigma\in \texttt{Aut}(R)$ is an \textbf{involution} if $\sigma^{2}=id_{R}$, or equivalently, if $\sigma=\sigma^{-1}$. For example, $R$ always admit the \textit{trivial involution} $id_{R}:x\longmapsto x$ and the \textit{involution sign change} $\sigma:x\longmapsto -x$. For each involution $\sigma\in \texttt{Aut}(R)$, consider the correspondence $\tilde{\sigma}:\texttt{Rd}(R)\longrightarrow \texttt{Rd}(R)$ given by $\tilde{\sigma}(I):=\sigma(I)$ for every $I\in \texttt{Rd}(R)$.

For each involution $\sigma\in \texttt{Aut}(R)$, $\tilde{\sigma}$ is well defined. Indeed, let $I\in \texttt{Rd}(R)$. Given $x\in \sqrt{\sigma(I)}$, we have $x^n\in \sigma(I)$ for some $n>0$. Thus, we get $x^{n}=\sigma(i)$ for some $i\in I$. Therefore, we have
$$
(\sigma(x))^{n}=\sigma(x^{n})=\sigma(\sigma(i))=\sigma^{2}(i)=i\in I
$$
and $\sigma(x)\in \sqrt{I}$. Since $I$ is a radical ideal, then $\sigma(x)\in I$. So, there is $i'\in I$ such that $\sigma(x)=i'$. Therefore, we have
$$
x=\sigma^{-1}(i')=\sigma(i')\in \sigma(I).
$$
This shows that $\sqrt{\sigma(I)}\subseteq\sigma(I)$ and, then, $\sigma(I)$ is a radical ideal.

Hence, if $G\subset \texttt{Aut}(R)$ is any group of involutions of $R$, then $\texttt{Rd}(R)$ is a $G$-invariant set. Therefore,
$H_{G}(\texttt{Rd}(R))$ is a $\omega$-stable group over $\texttt{Rd}(R)$.
\end{ex}

\begin{ex}
    Let $D$ be an integral domain, but not a field, and $R:=D[x]$ be the ring of polynomials in the variable $x$, with coefficients in $D$. For each $d\in D$, let be the ideal $I_{d}:=\left\langle x-d \right \rangle$ and consider the collection $\mathcal{I}_{d}:=\lbrace I_{d '}\,|\,d'\,\,\textrm{is associated to}\,\,d\rbrace$. Note that $\Omega(R,\mathcal{I})=\lbrace 1\rbrace$. In fact, since $\dfrac{R}{I_d}\simeq D$ and $D$ is an integral domain, then each $I_d$ is a prime ideal in $R$. Thus, each collection $\mathcal{I}_d$ only has a single $\omega_R$-orbit.
    
    Now, for each $\mu\in D^{\times}$, let us consider the automorphism
    $$
    g_{\mu}:d_{n}x^{n}+d_{n-1}x^{n-1}+\cdots +d_{1}x+d_0\longmapsto d_{n}x^{n }+d_{n-1}x^{n-1}+\cdots +d_{1}x+\mu d_0
    $$
    Note that each set $G:=\lbrace g_{\mu}\,|\,\mu\in D^{\times}\rbrace$ is a subgroup of $Aut(R)$. Furthermore, each collection $\mathcal{I}_{d}$ is $G$-invariant; in fact, let $I_{d'}\in \mathcal{I}_d$ and $\mu\in D^{\times}$. So, $d'=\theta d$, for some $\theta\in D^{\times}$, and we have
    
$$
g_{\mu}(I_{d'})=g_{\mu}(\left\langle x-d' \right \rangle)=\lbrace g_{\mu}(p(x)(x-d')) \,|\,p(x)\in R\rbrace=
$$
 
 $$
 =\lbrace g_{\mu}(p(x))g_{\mu}(x-d')\,|\,p(x)\in R\rbrace=\lbrace q(x)(x-\mu d')\,|\,q(x)\in R\rbrace\overset{d'=\theta d}{=}
 $$

$$
 \overset{d'=\theta d}{=}\lbrace q(x)(x-d'')\,|\,q(x)\in R\rbrace=\left \langle x-d" \right \rangle=I_{d"}\in \mathcal{I}_d.
 $$
 where $d"=\mu\theta d$.
 
 In this way, for each group $H_{G}(\mathcal{I}_d)$, every $H_{G}(\mathcal{I}_d)$-orbit is contained in the only $\omega$-orbit $[I_{d}]_{\omega}$ and, hence, each group $H_{G}(\mathcal{I}_d)$ is $\omega$-stable over $\mathcal{I}_{d}$.
 \end{ex}
Next, we provide a example of groups that are not $\omega_R$-stable.
\begin{ex}
   Let $R$ be an artinian ring, but not a field. We know that $\Omega(R)=\lbrace 1,2,\ldots,n\rbrace$, for some $n\in\mathbb{N}$, $n>1$. Making use of the \textit{ZFC-Axiom}, if necessary, we can take a non-empty set $\mathcal{I}$ of absorbing ideals of $R$, whose $\omega$-orbits are singletons. Thus, no non-trivial subgroup $H$ of $\texttt{Perm}(\mathcal{I})$ is $\omega$-stable over $\mathcal{I}$. More generally, we can take sets of absorbing ideals like these, via the ZFC-axiom, in any ring $R$ such that $\infty\not\in \Omega(R)$.
   \end{ex}

Let $X$ be a non-empty set. If $G$ is a group acting over $X$, we say that $X$ is a $G$-\textbf{set}. For each $g\in G$ and $x\in X$, let is denote by $g\cdot x$ the image of the action of $g$ on $x$ whenever we want to omit it or when there are no ambiguities about the action in question. 

 Let $\equiv$ be an equivalence relation on $X$. We say that $\equiv$ is a $G$-\textbf{congruence} if, for any $x,x'\in X$ and $g\in G$,
$$
x\equiv x'\Leftrightarrow g\cdot x\equiv g\cdot x'
$$
 It is not hard to see that the trivial equivalence relations given by 
 $$
 x\equiv_1 x'\Leftrightarrow x,x'\in X\,\,\,\,\textrm{and}\,\,\,\, x\equiv_2 x'\Leftrightarrow x=x'
 $$ 
 are the \textbf{trivial} $G$-\textbf{congruence}.

 Let $\mathcal{I}$ be a non-empty set of proper ideals of $R$ and $H$ be a subgroup of $\texttt{Perm}(\mathcal{I})$.
    
\begin{prop}\label{prop3.1}
    If $H$ is $\omega$-stable over $\mathcal{I}$ then $\sim_{\omega}$ is a $H$-congruence. 
\begin{proof}
  
Let $I,J\in\mathcal{I}$. Then 

$$
I\sim_{\omega} J \Leftrightarrow \omega(I)=\omega(J)
$$
If $H$ is $\omega$-stable then, for every $h\in H$, we have
$$
\omega(h(I))=\omega(I)\,\, \textrm{and}\,\, \omega(h(J))=\omega(J).
$$
Thus $\sim_{\omega}$ is a $H$-congruence. 
\end{proof}

\end{prop}
With the condition that $H$ is a $\omega$-stable group over $\mathcal{I}$, we characterize all $\sim_{\omega}$ equivalence relations that are non-trivial $H$-congruence. It is a direct consequence of the Proposition \ref{prop3.1}.
\begin{cor}\label{cor3.1}
    Let $H$ be a $\omega$-stable group over $\mathcal{I}$. Then $\sim_{\omega}$ is a non-trivial $H$-congruence iff the $\omega$-orbits of $\mathcal{I}$ are not trivial.
\end{cor}

Let $\equiv$ be an equivalence relation over $\mathcal{I}$. We say that $\equiv$ is a $\omega$-\textbf{congruence} if

$$
I\equiv J \Rightarrow I\sim_{\omega} J
$$
for every $I,J\in \mathcal{I}$.

\begin{rem}\label{obs3.3}
    Hence, it is clear that a subgroup $H$ of $\texttt{Perm}(\mathcal{I})$ is $\omega$-stable over $\mathcal{I}$ iff the natural equivalence relation $\equiv_{H}$ given by $H$-orbits is a $\omega$-congruence.
\end{rem}

Therefore, applying the Proposition \ref{prop3.1} and the Remark \ref{obs3.3} we have the follow result. 

\begin{cor}\label{cor3.2}
    If $\equiv_{H}$ is a $\omega$-congruence then $\sim_{\omega}$ is a $H$-congruence.
\end{cor}

\begin{quest}
    The Corollary \ref{cor3.2} given a necessary condition for a subgroup $H$ of $\texttt{Perm}(\mathcal{I})$ to be a $\omega$-stable group over $\mathcal{I}$. Now, for what set of proper ideals $\mathcal{I}$ and subgroups $H$ of $\texttt{Perm}(\mathcal{I})$ is worth the reciprocal of Corollary \ref{cor3.2}?
\end{quest}

\section{Transitive actions of $\omega$-stable groups}

We say that $G$ is \textbf{transitive} on $X$ if, for any $x,x'\in X$, there exist $g\in G$ such that $x'=g\cdot x$. Equivalently, $G$ is transitive over $X$ iff $X$ has only one $G$-orbit, i.e, $X=\texttt{Orb}_{G}(x)$, for some (and hence all) $x\in X$.

Let $\mathcal{I}$ be a set of proper ideals of $R$ and $H$ be a subgroup of $\texttt{Perm}(\mathcal{I})$.

\begin{prop}\label{prop.4.1}
    If $H$ is $\omega$-stable and transitive over $\mathcal{I}$, then  $\Omega(R,\mathcal{I})=\lbrace \omega(I)\rbrace$ for some ideal $I\in \mathcal{I}$.
   \begin{proof}
     If $H$ is $\omega$-stable on $\mathcal{I}$ then every $H$-orbit is contained in only one $\omega$-orbit of $\mathcal{I}$. Hence, since $H$ is transitive over $\mathcal{I}$, i.e, $\mathcal{I}=\texttt{Orb}_{H}(I)$, for some $I\in \mathcal{I}$, then by the Remark \ref{obs2.2} the result it is proven.
   \end{proof}
   \end{prop}
\begin{cor}\label{cor4.1}
 Let $R$ be an integral domain and suppose that $\mathbf{0}\in \mathcal{I}$. If $H$ is $\omega$-stable and transitive over $\mathcal{I}$, then all ideals in $\mathcal{I}$ are prime.
 \begin{proof}
     Suppose that $H$ is $\omega$-stable and transitive over $\mathcal{I}$. By Proposition \ref{prop.4.1}, $\Omega(R,\mathcal{I})=\lbrace \omega_{R}(I)\rbrace$ for some (and all) ideal $I\in \mathcal{I}$. Now, just choose $I=\mathbf{0}$.
 \end{proof}
   \end{cor}

\begin{cor}\label{cor4.2}
    Let $R$ be an integral domain and $G$ be a subgroup of involutions of $R$. If $H_{G}(\texttt{Rd}(R))$ is transitive over $\texttt{Rd}(R)$ then every radical ideal of $R$ is prime (primary).
    \begin{proof}
  Recall that a proper ideal of $R$ is primary if and only if it is radical is prime. By the Example \ref{ex3.1}, $H_{G}(\texttt{Rd}(R))$ is $\omega$-stable on $\texttt{Rd}(R)$. Since $\mathbf{0}=\sqrt{\mathbf{0}}$, so it follows from the Corollary \ref{cor4.1}.
    \end{proof}
   \end{cor}

   The next result tells us that the search for $\omega$-stable groups over the set of all proper ideals of integral domains that are not fields can and must be done by eliminating the transitivity condition.

\begin{cor}\label{cor4.3}
    Let $R$ be an integral domain. Then $H$ is $\omega$-stable and transitive over $\texttt{Prp}(R)$ iff $R$ is a field.
    \begin{proof}
        Suppose that $H$ is $\omega$-stable and transitive over $\texttt{Prp}(R)$. Since $\mathbf{0}\in \texttt{Prp}(R)$ so, by Corollary \ref{cor4.1}, all ideals of $R$ are prime. Being all proper ideals of $R$ primes, let any $x\in R$, $x\neq 0$. Since $R$ is an integral domain, then $x^{2}\neq 0$. Take the non-zero ideal $\left\langle x^{2}\right\rangle$. Since $\left\langle x^{2}\right\rangle$ is prime (by hypothesis) and $x^{2}\in \left\langle x^{2}\right\rangle$, then $x\in \left\langle x^{2}\right\rangle$. Hence, $x=x^{2}y$, for some $y\in R$. Thus, we have
        $$
        x(xy-1)=0
        $$
        Since $R$ is an integral domain and $x\neq 0$, we obtain $xy=1$. Thus, $x$ is invertible. Therefore, $R$ is a field. The reciprocal is obvious.
    \end{proof}
\end{cor}

\centering

\end{document}